\documentclass[leqno]{article}

\begin{document}

\title{Potpourri, 8}

\author{Stephen William Semmes	\\
	Rice University		\\
	Houston, Texas}

\date{}

\maketitle

	Let $E$ be a finite nonempty set, and let $\mathcal{F}(E, {\bf
R})$, $\mathcal{F}(E, {\bf C})$ denote the vector spaces of real and
complex-valued functions on $E$.  Here we shall be primarily
interested in the norms
\begin{equation}
	\|f\|_1 = \sum_{x \in E} |f(x)|
\end{equation}
and
\begin{equation}
	\|f\|_\infty = \max \{|f(x)| : x \in E\}
\end{equation}
for functions $f$ on $E$.

	Let $V$ be a finite-dimensional real or complex vector space.
A linear transformation $A$ from real or complex-valued functions on
$E$ into $V$, as appropriate, is characterized by the images under $A$
of the functions on $E$ which are equal to $1$ at one element of $E$
and equal to $0$ elsewhere on $E$.  A linear transformation $T$ from
$V$ into real or complex-valued functions on $E$, as appropriate, is
basically the same as a collection of linear functionals on $V$, one
for each element of $E$, corresponding to evaluating images of vectors
in $V$ under $T$ at elements of $E$.

	Suppose that $\|\cdot \|_V$ is a norm on $V$.  Thus $\|v\|_V$
is a nonnegative real number for all $v \in V$, $\|v\|_V = 0$ if and
only if $v = 0$, $\|\alpha \, v \|_V = |\alpha| \, \|v\|_V$ for all
real or complex numbers $\alpha$, as appropriate, and all $v \in V$,
and
\begin{equation}
	\|v + w\|_V \le \|v\|_V + \|w\|_V
\end{equation}
for all $v, w \in V$.

	Let $A$ be a linear transformation from functions on $E$ into
$V$.  For each $x \in E$, let $\delta_x$ be the function on $E$ given
by $\delta_x(x) = 1$, $\delta_x(y) = 0$ when $y \ne x$.  It is easy to
see that
\begin{equation}
	\|A(f)\|_V \le \max \{\|A(\delta_x)\|_V : x \in E \} \, \|f\|_1
\end{equation}
for all real or complex-valued functions $f$ on $E$, as appropriate.

	If $\lambda$ is a linear functional on $V$, which is to say a
linear mapping from $V$ into the real or complex numbers, as
appropriate, then we put
\begin{equation}
	\|\lambda\|_{V^*} = \sup \{|\lambda(v)| : v \in V, \ \|v\|_V \le 1 \}.
\end{equation}
It is well-known that this defines a norm on $V^*$, the vector space
of linear functionals on $V$.

	If $T$ is a linear mapping from $V$ into functions on $E$,
then for each $x \in E$ we get a linear functional $\lambda_x$ on $V$,
given by
\begin{equation}
	\lambda_x(v) = T(v)(x)
\end{equation}
for all $x \in E$ and $v \in V$.  It is easy to see that
\begin{equation}
	\|T(v)\|_\infty \le \max \{\|\lambda_x\|_{V^*} : x \in E\} \, \|v\|_V
\end{equation}
for all $v \in V$.

	A famous result states that if $W$ is a linear subspace of $V$,
$\lambda$ is a linear mapping from $W$ into the real or complex numbers,
as appropriate, and $L$ is a nonnegative real number such that
\begin{equation}
	|\lambda(w)| \le L \, \|w\|_V
\end{equation}
for all $w \in W$, then there is an extension of $\lambda$ to a linear
functional on all of $V$ which satisfies the same inequality.  This
extension result also works for linear mappings into functions on $E$,
using the norm $\|f\|_\infty$ for functions on $E$.  Indeed, one can
apply the result for linear functionals to each component of such a
mapping.

	Let $\Sigma$, $\Sigma^*$ denote the unit spheres in $V$, $V^*$
with respect to the norms $\|v\|_V$, $\|\lambda\|_{V^*}$, i.e.,
\begin{equation}
	\Sigma = \{v \in V : \|v\|_V = 1\}
\end{equation}
and
\begin{equation}
	\Sigma^* = \{\lambda \in V^* : \|\lambda\|_{V^*} = 1\}.
\end{equation}
Let $C(\Sigma)$, $C(\Sigma^*)$ denote the vector spaces of real or
complex-valued continuous functions on $\Sigma$, $\Sigma^*$, according
to whether $V$ is real or complex.

	For each $v \in V$, let $\phi_v$ be the continuous function on
$\Sigma^*$ such that $\phi_v(\lambda) = \lambda(v)$.  Thus $v \mapsto
\phi_v$ defines a linear mapping of $V$ into $C(\Sigma^*)$.  We also
have that
\begin{equation}
	\max \{|\phi_v(\lambda)| : \lambda \in \Sigma^*\} = \|v\|_V,
\end{equation}
because for each $v \in V$ there is a $\lambda \in \Sigma^*$ such that
$\lambda(v) = \|v\|_V$.

	Let us write $M(\Sigma)$ for the space of \emph{measures} on
$\Sigma$, which means the space of linear mappings $\mu$ from
$C(\Sigma)$ into the real or complex numbers, according to whether $V$
is real or complex, which are bounded in the sense that there is a
nonnegative real number $k$ such that
\begin{equation}
	|\mu(f)| \le k \, \max \{|f(v)| : v \in \Sigma\}
\end{equation}
for all $f \in C(\Sigma)$.  In this event we put
\begin{equation}
	\|\mu\|_M = \sup \{|\mu(f)| : f \in C(\Sigma), \ 
				|f(v)| \le 1 \hbox{ for all } v \in \Sigma\},
\end{equation}
which defines a norm on $M(\Sigma)$, and which is the same as the
smallest possible value of $k$ in the preceding inequality.

	A measure $\mu$ on $\Sigma$ is said to be nonnegative if
$\mu(f)$ is a nonnegative real number whenever $f$ is a real-valued
continuous function on $\Sigma$ such that $f(v) \ge 0$ for all $v \in
\Sigma$.  For a nonnegative measure the boundedness condition is
automatic, and $\|\mu\|_M$ is equal to $\mu$ applied to the constant
function equal to $1$ on $\Sigma$.

	If $\mu$ is a measure on $\Sigma$ and $\phi$ is a continuous
function on $\Sigma$, then we get a new measure $\mu_\phi$ defined by
\begin{equation}
	\mu_\phi(f) = \mu (\phi \, f)
\end{equation}
for all $f \in C(\Sigma)$.  Notice that
\begin{equation}
	\|\mu_\phi\|_M \le \max \{|\phi(v)| : v \in \Sigma\} \, \|\mu\|_M.
\end{equation}
If $\phi_1, \ldots, \phi_l$ are continuous functions on $\Sigma$, then
\begin{equation}
	\sum_{j=1}^l \|\mu_{\phi_j}\|_M 
		\le \max \bigg\{ \sum_{j=1}^l |\phi_j(v)| : v \in \Sigma\bigg\}
			\, \|\mu\|_M.
\end{equation}

	Let $\mu$ be a measure on $\Sigma$.  We can extend $\mu$ to a
linear mapping from $V$-valued continuous functions on $\Sigma$ into
$V$ in a canonical way, so that $\mu$ applied to the product of a
continuous scalar-valued function $f$ on $\Sigma$ and a vector $v \in
V$ is equal to $\mu(f)$ times $v$.

	Suppose that $f$ is a continuous $V$-valued function on
$\Sigma$ and that $\lambda$ is a linear functional on $\Sigma$.  We
can also characterize $\mu(f) \in V$ by saying that $\lambda(\mu(f))$
is equal to $\mu$ applied to the scalar valued function
$\lambda(f(v))$ on $\Sigma$.  Using this characterization it follows
easily that
\begin{equation}
	\|\mu(f)\|_V \le \max \{\|f(v)\|_V : v \in \Sigma\} \, \|\mu\|_M
\end{equation}
for all continuous $V$-valued functions $f$ on $\Sigma$.

	Define a linear mapping $\Phi$ from $M(\Sigma)$ into $V$ by
saying that $\Phi$ applied to a measure $\mu$ on $\Sigma$ is equal to
$\mu$ applied to the $V$-valued function on $\Sigma$ which sends each
$v \in \Sigma$ to itself.
Thus we have that
\begin{equation}
	\|\Phi(\mu)\|_V \le \|\mu\|_M
\end{equation}
for every measure $\mu$ on $\Sigma$.

	For each $v \in \Sigma$ we can define the measure on $\Sigma$
which is the Dirac mass at $v$, and which sends a function on $\Sigma$
to its value at $v$.  The Dirac mass at $v$ has norm $1$, and $\Phi$
applied to the Dirac mass at $v$ is equal to $v$.

	Of course $C(\Sigma)$ is a Banach space with respect to the
supremum norm
\begin{equation}
	\|f\| = \max \{|f(v)| : v \in \Sigma\}.
\end{equation}
This space is separable, meaning that it has a countable dense subset.
This can be shown using the fact that $\Sigma$ is compact and
continuous functions on $\Sigma$ are uniformly continuous.

	Let us recall that a closed and bounded subset of $C(\Sigma)$
is compact if and only if it is equicontinuous.  This is an instance
of the theorem of Arzela and Ascoli.

	The space $M(\Sigma)$ of measures on $\Sigma$ is a Banach
space with respect to the norm $\|\mu\|_M$.  Alternatively, it is
frequently natural to use another topology on $M(\Sigma)$, which is
the weakest topology in which
\begin{equation}
	\mu \mapsto \mu(f)
\end{equation}
is a continuous function on $M(\Sigma)$ for each $f \in C(\Sigma)$.
With respect to this topology, $M(\Sigma)$ is still a nice
locally convex topological vector space.

	One can restrict this weak topology to the closed unit ball in
$M(\Sigma)$ defined by the norm, or to any other closed ball for that
matter.  Although the topology is not defined by a single norm, its
restriction to a closed ball is metrizable, and the closed ball
becomes compact in this topology, by well-known results.


\begin{thebibliography}{14}


\bibitem {F} G.~Folland, {\it Real Analysis: Modern Techniques and
their Applications}, 2nd edition, Wiley, 1999.

\bibitem {G} R.~Goldberg, {\it Methods of Real Analysis}, 2nd edition,
Wiley, 1976.

\bibitem {H-S} E.~Hewitt and K.~Stromberg, {\it Real and Abstract
Analysis}, Springer-Verlag, 1975.

\bibitem {Frank} F.~Jones, {\it Lebesgue Integration on Euclidean
Space}, Jones and Bartlett, 1993.

\bibitem {K} S.~Krantz, {\it Real Analysis and Foundations},
CRC Press, 1991.

\bibitem {L-T1} J.~Lindenstrauss and L.~Tzafriri, {\it Classical
Banach Spaces I, Sequence Spaces}, Springer-Verlag, 1977.

\bibitem {L-T2} J.~Lindenstrauss and L.~Tzafriri, {\it Classical
Banach Spaces II, Function Spaces}, Springer-Verlag, 1979.

\bibitem {Re} M.~Reed, {\it Fundamental Ideas of Analysis},
Wiley, 1998.

\bibitem {Ro} H.~Royden, {\it Real Analysis}, 3rd edition, MacMillan,
1988.

\bibitem {Ru1} W.~Rudin, {\it Principles of Mathematical Analysis},
3rd edition, McGraw-Hill, 1976.

\bibitem {Ru2} W.~Rudin, {\it Real and Complex Analysis}, 3rd edition,
McGraw-Hill, 1987.

\bibitem {Ru3} W.~Rudin, {\it Functional Analysis}, 2nd edition,
McGraw-Hill, 1991.

\bibitem {KS} K.~Stromberg, {\it Introduction to Classical
Real Analysis}, Wadsworth, 1981.

\bibitem {RS} R.~Strichartz, {\it The Way of Analysis}, revised edition,
Jones and Bartlett, 2000.

\end{thebibliography}
\end{document}